\documentclass[11pt]{article}

 \usepackage{amsfonts,amssymb,graphicx}

\usepackage[colorlinks=true,linkcolor=blue,citecolor=blue]{hyperref}

\newtheorem{rema}{Remark}
\newtheorem{defi}{Definition}
\newtheorem{lemm}{Lemma}
\newtheorem{theo}{Theorem}
\newtheorem{coro}{Corollary}
\newtheorem{exem}{Example}

\newcommand{\Z}[1][]{\ensuremath{{\mathbb{Z}^{#1}} }}
\newcommand{\C}[1][]{\ensuremath{{\mathbb{C}^{#1}} }}
\newcommand{\R}[1][]{\ensuremath{{\mathbb{R}^{#1}} }}
\renewcommand{\S}[1][]{\ensuremath{{\mathbb{S}^{#1}} }}

\newcommand{\X}[1][]{\ensuremath{{\mathbb{X}^{#1}} }}

\def\Re{ \mathrm{Re}\, }
\def\Im{ \mathrm{Im}\, }

\renewcommand{\j}{\ensuremath{J}}
\renewcommand{\L}{{\cal L}}

\newcommand{\s}{{\cal S}}
\newcommand{\<}{\langle}
\renewcommand{\>}{\rangle}
\newcommand{\ga}{\gamma}

\newcommand{\eps}{\epsilon}

\newcommand{\la}{\lambda}
\newcommand{\be}{\beta}

\date{}

\title{Minimal  Lagrangian submanifolds in indefinite complex space}
\author{ Henri Anciaux\footnote{The author is supported by CNPq (PQ 302584/2007-2)}}
\begin{document}
\maketitle

\centerline{\textbf {Abstract}}

\medskip

{\small \noindent Consider the  complex linear space $\Bbb C^n$ endowed with the canonical 
pseudo-Hermitian form of signature $(2p,2(n-p)),$
where $0 \leq p \leq n.$
This yields both a pseudo-Riemannian and a symplectic structure on
$\C^n$. We prove that those submanifolds which are both Lagrangian
and minimal with respect to these structures minimize the volume in
their Lagrangian homology class. We also describe several families of
minimal Lagrangian submanifolds. In particular, we characterize
the minimal Lagrangian surfaces in $\Bbb C^2$ endowed with its natural neutral metric and the
equivariant minimal Lagrangian submanifolds of $\Bbb C^n$ with arbitrary signature.}

\medskip

\centerline{\small \em 2000 MSC: 53D12,  49Q05
\em }

\section*{Introduction}

It has been discovered in the seminal paper of Harvey and Lawson
\cite{HL1} (see also \cite{Ha}) that a minimal Lagrangian
submanifold of complex Euclidean space is calibrated and therefore
minimizes the area in its homology class. This remarkable fact no
longer holds true in an arbitrary K\"ahler manifold, but it does
in a certain class of K\"ahler manifolds, namely the \em
Calabi-Yau manifolds. \em The  study of minimal Lagrangian
submanifolds (usually called \em Special Lagrangian submanifolds\em)
 in Calabi-Yau manifolds has attracted much attention
recently, in particular because of its close relationship with
mirror symmetry, a important issue in theoretical physics
(\cite{SYZ}).

\medskip

Most of the theory of submanifolds in Riemannian geometry may be
extended to the realm of pseudo-Riemannian geometry, and recently
this issue has attracted much attention. In particular, Mealy extended in
\cite{Me} (see also \cite{HL2}) the concept of calibration in pseudo-Riemannian
manifolds. On the other hand, Dong
addressed in \cite{Do} the local study of minimal Lagrangian submanifolds in
complex linear space $\C^n$ endowed with the pseudo-Hermitian form
defined by:
$$\<\<.,.\>\>_p := -\sum_{j=1}^{p}  dz_j d\bar{z}_j+\sum_{j=p+1}^{n}  dz_j d\bar{z}_j,$$
where $ 0 \leq p \leq n.$ If $p=0$ or $n$, we fall back in the
classical, Riemannian setting of \cite{HL1}.
 One of Dong's main observations is
that, although the geometry of a minimal Lagrangian in
$(\C^n,\<\<.,.\>\>_p), p \neq 0,n$ is somehow analogous to that of the Riemannian
case $p=0,$ they are always unstable (in the classical sense), so
in particular they can not be homology minimizing.

\medskip

The main result of this paper is
 that although the original calibration of Harvey and
Lawson does not calibrate minimal Lagrangian submanifolds (as
pointed by Dong), it does calibrate them in their Lagrangian
homology class (Main Theorem, Section \ref{secmini}). 
 In the remainder of the paper we describe some
families of minimal Lagrangian submanifolds. In particular, we
show that a minimal Lagrangian surface of $(\C^2,\<\<.,.\>\>_1)$
must be the Cartesian product of two curves contained in two
mutually orthogonal, null, non-Lagrangian planes (Theorem
\ref{minlagn=2}, Section \ref{plane}). We also characterize $SO(p,n-p)$-equivariant
minimal Lagrangian submanifolds of $(\C^n,\<\<.,.\>\>_p)$
(Section \ref{secequi}). This family generalizes the \em Lagrangian
catenoid, \em which was first described by Harvey and Lawson and
 studied in more detail in \cite{CU}. Finally, inspired by a
construction due to Joyce (\cite{Jo}), we produce a larger family
of minimal Lagrangian submanifolds obtained from evolving quadrics
(Section \ref{secquad}).

\section{Preliminaries}
Consider the complex linear space $\C^n$ of arbitrary dimension
$n$, endowed with its canonical complex structure $J$ and, for $0 \leq p \leq n,$
 the pseudo-Hermitian form of arbitrary signature $(p,n-p)$ defined by:
$$\<\<.,.\>\>_p := -\sum_{j=1}^{p}  dz_j d\bar{z}_j+\sum_{j=p+1}^{n}  dz_j d\bar{z}_j.$$
The real and imaginary parts of $\<\<.,.\>\>_p$
$$ \<.,.\>_{2p}=\Re \<\<.,.\>\>_p \quad \mbox{ and } \quad  \omega_p= -\Im \<\<.,.\>\>_p$$
 yield two different
structures: while the bilinear form $ \<.,.\>_{2p}$ is a
pseudo-Riemannian metric with signature $(2p,2(n-p)),$ the closed
$2$-form $\omega_p$ is, up to a (real) linear change of
coordinates, the canonical symplectic form of $\C^n \simeq T^*
\R^n$ regarded as the cotangent bundle of $\R^n.$

A smooth immersed submanifold ${\s}$ of $(\C^n,\<\<.,.\>\>_p)$  is
said to be \em non-degenerate \em if the induced metric on $\s$ is
itself non-degenerate. Moreover, a non-degenerate submanifold is
said to be \em minimal \em if it is a critical point of the  volume  with respect
to compactly supported variations. On the other hand, an
$n$-dimensional submanifold is said to be \em Lagrangian \em if
$\omega_p$ vanishes on it. The equation $\omega_p = \<J.,.\>_{2p}$
shows that a non-degenerate submanifold is Lagrangian if and only if its tangent
and normal bundles $T\L$ and $N \L$ are isometrically exchanged by the complex
structure $J.$ Since $ T\L \oplus N \L = T \C^n$, the following fact holds:
\begin{lemm} \label{signature} The induced metric on a non-degenerate Lagrangian
submanifold of $(\C^n,\<\<.,.\>\>_p)$ has signature $(p,n-p).$
\end{lemm}

 We furthermore introduce the
   holomorphic volume form $\Omega := dz_1 \wedge ... \wedge dz_n,$
   which turns out to be useful for the description of the
   geometry of a Lagrangian submanifold.

\begin{defi}  The \em Lagrangian angle\index{Lagrangian angle} \em $\be$
of a non-degenerate, Lagrangian, oriented,  submanifold $\L$  is the
map $\be: \L \to \R/ 2\pi \Z$ defined by
$$ \be:=\arg \Omega(X_1,...,X_n),$$
where $(X_1,...,X_n)$ is a tangent moving frame along $\L$ (it is
easy to check that the definition of $\be$ does not
depend on the choice of the moving frame, see \cite{An}).
\end{defi}

The importance of the Lagrangian angle map is due to the following
formula, which was first derived by Chen and Morvan  in the definite case (see
\cite{CM}) and extended to the indefinite case in
(\cite{Do}):

\begin{theo} \label{Hlag} Let $\L$ be a non-degenerate, Lagrangian submanifold of $(\C^{n},\<\<.,.\>\>_p)$
with Lagrangian angle $\be$ and mean curvature vector $\vec{H}$.
Then the following formula holds
$$ n\vec{H} = J \nabla \be,$$
where $\nabla$ denotes the gradient operator with respect to  the induced metric.
\end{theo}

\begin{coro} \label{corolag}
A Lagrangian submanifold $\L$ of $(\C^n,\<\<.,.\>\>_p)$ is minimal if and only if
it has constant Lagrangian angle.
\end{coro}

\section{Minimizing properties}  \label{secmini}
We recall here the simple but powerful concept developed in \cite{HL1}: let $({\cal M},g)$ be a Riemannian manifold. A \em calibration \em is a closed
 $n$-form $\Theta$ of ${\cal M}$ which is bounded from above by the $n$-th
 dimensional volume form induced from $g,$ i.e., for any
 $n$-vector $X_1 \wedge ... \wedge X_n$, we have
 $$\Theta(X_1 , ..., X_n) \leq \sqrt{\left|\mbox{det}_{\R}  [g(X_j,X_k)]_{1 \leq j,k \leq n}\right|}:= d\mbox{Vol}(X_1,...,X_n) .$$
  A $n$-dimensional
 submanifold $\s$ of ${\cal M}$ is said to be \em calibrated \em by $\Theta$
 if the restriction of $\Theta$ to $\s$ is equal to the
 $n$-volume, i.e.\ equality is attained in the expression above when
 $(X_1,...,X_n)$ is a tangent moving frame along $\s.$ By Stokes
 theorem, it follows that if $\s'$ is any submanifold
 belonging to the homology class of a calibrated submanifold
 $\s$, we have
 $$ \mbox{ Vol}(\s) = \int_{\s} d \mbox{Vol} =  \int_{\s}
 \Theta = \int_{\s'} \Theta \leq   \int_{\s'} d \mbox{Vol}
 = \mbox{ Vol}(\s').$$
 Therefore a calibrated submanifold minimizes the volume in its homology class, hence
 it is in particular minimal and stable.

   \medskip

Among the few known examples of calibrations is the $1$-parameter family of $n$-forms of $\C^n,$
            $$\Theta_0:=\Re (e^{-i \be_0} \Omega), \, \, \be_0 \in \R / 2 \pi \Z,$$
             discovered in \cite{HL1}. The
           calibrated submanifolds of $\Theta_0$ are precisely those Lagrangian
 submanifold of   $(\C^n, \<\<.,.\>\>_0)$ with constant Lagrangian angle $\be_0.$

On the other hand, it was proved in \cite{Do} (see also \cite{An}) that
a minimal submanifold of  $(\C^n, \<\<.,.\>\>_p)$ whose tangent or  normal bundle is indefinite is unstable.
 By Lemma \ref{signature}, it follows that
there is no hope to find a calibration for minimal Lagrangian submanifolds of  $(\C^n, \<\<.,.\>\>_p)$ in the usual sense when $p \neq 0,n.$

Nevertheless,  the following result holds:

\bigskip

\noindent \textbf{Main Theorem:} \em Let $\L$ be a minimal Lagrangian submanifold of $(\C^n, \<\<.,.\>\>_p).$ Then $\L$ minimizes
the volume in its Lagrangian homology class. \em

\bigskip

\noindent \textit{Proof.}
Let $\be_0$ be the constant Lagrangian angle of $\L.$
We claim that if $X_1, ...,X_n$ are $n$ vectors spanning a non-degenerate Lagrangian subspace, then
  $ \Theta_0(X_1, ..., X_n)  \leq d\mbox{Vol}(X_1,...,X_n),$
  with equality if and only if $\beta(X_1,...,X_n)=\be_0.$
  To see this, observe that given a vector $X$ of $\C^n,$ we have
$$ X = \sum_{j=1}^n \eps_j \<\<X,e_j\>\>_p e_j,$$
where  $(e_1,...,e_n)$ is the canonical Hermitian basis of $(\C^n,\<\<.,.\>\>_p)$ and $\eps_j:=\<\<e_j,e_j\>\>_p=\pm 1.$ 
Setting $M:= [ \<\<X_j, e_k\>\>_p]_{1 \leq j,k \leq n},$ it follows that
$$|\Omega (X_1, ...,X_n)|=|\mbox{det}_{\C} [\eps_k \<\<X_j, e_k\>\>_p]_{1 \leq j,k \leq n}|=|\mbox{det}_{\C} M|.$$
On the other hand, by the Lagrangian assumption, we have
\begin{eqnarray*}
 \<X_j ,X_k \>_{2p} &=&\<\<X_j ,X_k\>\>_p \\
&=& \sum_{l=1}^n   \<\<X_j, e_l\>\>_p\<\<X_k, e_l\>\>_p .
\end{eqnarray*}
Therefore
\begin{eqnarray*}
 d \mbox{Vol} (X_1, ...,X_n)&=& \left|\mbox{det}_{\R} ([\<X_j,X_k\>_{2p}]_{1 \leq j,k \leq n}) \right|^{1/2}\\
 &=& \left| \mbox{det}_{\R} (M . M^\ast) \right|^{1/2}\\
 &=& \left| \mbox{det}_{\C} (M . M^\ast) \right|^{1/2}\\
 &=& \big|\mbox{det}_{\C} M \big|,
 \end{eqnarray*}
where $M^\ast$ denotes the complex transpose of $M.$  It follows that
 $$ \Theta_0(X_1, ...,X_n) \leq \left|  \Omega(X_1,...,X_n) \right| =d \mbox{Vol}(X_1,...,X_n),$$
 and of course equality holds if and only if
 $$\be(X_1,...,X_n)=\arg (\Omega(X_1,...,X_n))=\be_0.$$
To conclude  the proof, we proceed exactly as in the case of a classical calibration: given a Lagrangian submanifold $\L'$ in the homology class of
$\L,$ we have
$$ \mbox{ Vol}(\L) = \int_{\L} d \mbox{Vol} =  \int_{\L}
 \Theta_0 = \int_{\L'} \Theta_0 \leq   \int_{\L'} d \mbox{Vol}
 = \mbox{ Vol}(\L').$$

\section{Examples of minimal Lagrangian surfaces in complex space}
 \subsection{Minimal Lagrangian surfaces in complex Lorentzian plane} \label{plane}
 In this section we characterize minimal Lagrangian surfaces of $\C^2$
 endowed with the "Lorentzian Hermitian metric"
  $$\<\<.,.\>\>_1 := -  dz_1 d\bar{z}_1+  dz_2 d\bar{z}_2=\<.,.\>_{2} - i \omega_1.  $$

\begin{theo} \label{minlagn=2} Let $\L$ be a minimal Lagrangian surface of $(\C^2,\<\<.,.\>\>_1).$ Then $\L$ is the product $\ga_1 \times \j \ga_2 \subset P \oplus \j P $, where  $\ga_1$
and $\ga_2$ are two planar curves contained in  a non-Lagrangian (and
therefore non-complex) null plane $P.$
\end{theo}

  \begin{lemm} \label{lemmn=2}
Let $P$ be a plane of $(\C^2, \<\<.,.\>\>_1).$ Then the induced metric on $P$ is totally null
  (i.e. $\<.,.\>_2|_{P}=0$)
  if and only if $\j P = P^{\omega_1}$, where
  $P^{\omega_1}$ denotes the symplectic orthogonal of $P.$

  \end{lemm}

  \noindent \textit{Proof.} Suppose first that $P$ is totally null and let $X$ be a vector  of $P$. For all vector $Y$ in $P$, we have
  $$0 =\<X,Y\>_2 = -\omega_1 (\j X,Y),$$
  so
  $\j X \in P^{\omega_1}$. Since it holds $\forall X \in P,$ we deduce that $\j P \subset P^{\omega_1}$, and
  the two-form $\omega_1$ being non-degenerate, $P^{\omega_1}$ is a two-dimensional subspace. Hence $\j P=P^{\omega_1}.$

Conversely, if $\j P= P^{\omega_1},$ then, for all  vector $X$ in $P$, we have
$ |X|^2_2 = -\omega_1(\j X,X)= 0.$ 
By the polarization formula
$2\<X,Y\>_2=|X+Y|_2^2 - |X|^2_2 - |Y|^2_2$, it implies that $P$ is totally null. 

\begin{rema} \em
This lemma proves in particular that a plane may be both complex and Lagrangian.
This fact may sound strange to the reader familiar with K\"{a}hler geometry,
 where complex and Lagrangian planes are two distinct classes.
 More precisely, if a plane enjoys any two of the three properties: \{totally null, Lagrangian, complex\}, then the third one holds as well.\em
\end{rema}

\noindent \textit{Proof of Theorem \ref{minlagn=2}.}

  Let $f: \L \to \C^2$ be a local parametrization of a minimal
  Lagrangian surface of $(\C^2, \<\<.,.\>\>_1).$
By  Lemma \ref{signature}, the induced metric on $\L$  is
Lorentzian, so it enjoys null coordinates $(u,v)$ (see \cite{We}).
A straightforward computation (see \cite{An}) shows that
$\vec{H}= \frac{2f_{uv}^\perp}{\<f_u,f_v\>_2},$ where $(.)^\perp$ denotes the projection onto the normal space. 
Moreover, differentiating the assumptions $|f_u|^2_2 = |f_v|_2^2=0,$ we get that that $f_{uv}$ is normal to $\L.$ 
It follows that the immersion $f$ is minimal if and only if $f_{uv}$ vanishes. Hence $f$ must take the form
$$f(u,v)= \ga_1(u) + \tilde{\ga}_2(v),$$ where $\ga_1,\tilde{\ga}_2$ are two  curves of $\C^2.$
Moreover, the assumption that $(u,v)$ are null coordinates translates into the fact that the two curves have null (i.e.\ lightlike)
velocity vector, and the non-degeneracy assumption is
$$\<\ga'_1(u),\tilde{\ga}'_2(v)\>_2 \neq 0, \forall \, (u,v) \in I_1 \times I_2.$$
On the other hand the Lagrangian assumption is:
$$\omega_1(\ga'_1(u),\tilde{\ga}'_2(v)) = 0, \forall \, (u,v) \in I_1 \times I_2.$$
The remainder of the proof relies on the analysis of the dimension of the two linear spaces $P_1 := Span \{ \ga'_1(u), u \in I_1 \}$ and $P_2 :=Span \{ \tilde{\ga}'_2(v), v \in I_2 \}.$ We first observe that $\dim P_1, \dim P_2 \geq 1$ and that the case $\dim P_1 = \dim P_2 =1$ corresponds to
the trivial case of $\L$ being planar.
Since the r\^{o}les of $\ga_1$ and $\tilde{\ga}_2$ are symmetric, we may assume without loss of generality that
$\dim P_1 \neq 1$.

Next, the Lagrangian assumption is equivalent to
 $P_2 \subset P_1^{\omega_1} $ and $P_1 \subset P_2^{\omega_1},$
so $\dim P_2 \leq \dim P_1^{\omega_1}$ and $\dim P_1 \leq \dim P_2^{\omega_1} $.
By the non-degeneracy of $\omega_1,$ it follows that
$\dim P_1 \leq \dim P_2^{\omega_1} = 4 - \dim P_2 \leq 3.$ We claim that in fact $\dim P_1= 2$. To see this,
assume by contradiction that $\dim P_1 = 3$. It follows that $\dim P_2 \leq \dim P_1^{\omega_1} =1,$ so
     the curve $\tilde{\ga}_2$ is a straight line, which may be parametrized as follows: $\tilde{\ga}_2(v)=e_0 v, $ where $e_0$ is a null vector of $\C^2.$ Then ${\ga}_1'$ is contained in the intersection of the light cone $\{(z_1,z_2) \in \C^2 \big| \, |z_1|=|z_2| \}$ with the hyperplane $\{e_0\}^{\omega_1}.$ An easy computation, using the fact that $e_0$ is null,
    shows that 
    $$\{(z_1,z_2) \in \C^2 \big| \, |z_1|=|z_2| \} \cap \{e_0\}^{\omega_1}=\Pi_1 \cup \Pi_2,$$ 
    where $\Pi_1$ and $\Pi_2$ are two null planes. Moreover, one of these  planes, say $\Pi_2$,   is contained in the metric orthogonal of $e_0$. By the non-degeneracy assumption
    $$\<\ga'_1(u),\tilde{\ga}'_2(v)\>_2 =
    v \<\ga'_1(u), e_0\>_2 \neq 0,$$
    we deduce that $\ga'_1 \in \Pi_1,$ which implies that $\dim P_1 \leq 2,$ a contradiction.

    To conclude, observe  that, using Lemma \ref{lemmn=2},
    $\tilde{\ga}_2 \in P_2 \subset P_1^{\omega_1} = \j P_1.$ Hence we just need to set $P:=P_1$ and $\ga_2:= - \j \tilde{\ga}_2,$ to get that
     $\ga_1, \ga_2 \subset P,$ so that $\L$ takes the required expression. 

 \subsection{Equivariant Lagrangian submanifolds in $\C^n$} \label{secequi}
\label{equivLagCn} In this subsection, we give a characterization
of those minimal Lagrangian submanifolds of $(\C^n,\<\<.,.\>\>_p)$
which are equivariant with respect to the canonical action of the
group $SO(p,n-p)$ defined as follows: for $z=x+iy \in \C^n$ and $M
\in SO(p,n-p)$ we  simply set $Mz := Mx + iMy.$ Of course we have
$\<\<Mz,Mz'\>\>_p=\<\<z,z'\>\>_p,$ so $SO(p,n-p)$ can be
identified with a subgroup of
$$U(n-p,p):=\{ M \in Gl(\C^{n}) | \,
\<\<MX,MY\>\>_p= \<\<X,Y\>\>_p \}.$$ Observe that  the orbits of
the action $SO(p,n-p)$ on $\R^n$ are the quadrics
$$ \X_{p,c}^{n-1}:=\big\{ x \in \R^n | \,  \<x,x\>_p=c \big\}.$$

\begin{theo}  \label{equivLagCnTheo}
Let $\L$ be an $SO(p,n-p)$-equivariant Lagrangian submanifold of
$\C^n$. Then it is locally congruent to the image of an immersion of the form
$$ \begin{array}{lccc}  f :
 & I \times \X^{n-1}_{p,\eps}
 &\to&  \C^n \\
& (s,x)& \mapsto &  \ga(s) x,
\end{array}$$
  where $\eps=1$ or $-1$ and $\ga: I \to \C^*$ is a planar curve. Moreover, the
Lagrangian angle of $\L$ is given by
$$ \be=\arg (\ga' {\ga}^{n-1}).$$
\end{theo}

\begin{rema} \em
In the definite, two-dimensional case ($p=0,n=2$), the $SO(2)$-action
 mentioned in the theorem above is not the only possible one, and there do
 exist Lagrangian surfaces of $\C^2$ equivariant by \em another \em $SO(2)$-action.
 For example, let $\ga(s)=(\ga_1(s),\ga_2(s))$ be a regular curve of
the sphere $\S^3$ such that  $\<\ga',J\ga\>_0 \neq 0$. Then the map $f(s,t) =(\ga_1(s)e^{it},\ga_2(s)e^{it})$ is a Lagrangian immersion
which is equivariant
 by the action $ M(z_1,z_2)=(Mz_1,Mz_2), M \in SO(2).$ These surfaces have been studied in \cite{Pi}, where
they are called \em Hopf surfaces.  \end{rema}

\noindent \textit{Proof of Theorem \ref{equivLagCnTheo}.}

\smallskip

\textbf{First case: $n =2$.}

\noindent Recall that the metric of $\C^2$ is
 $\<\<.,.\>\>_p=\eps_1 dz_1
d\bar{z}_1 + dz_2 d\bar{z}_2$ with $\eps_1=1$ or $-1.$ Introducing
$M_{\eps_1}:=\left(
\begin{array}{cc} 0 & -\eps_1
\\ 1 & 0
\end{array} \right),$ we have
$$ SO(2)= \{ e^{M_1 t}, t \in \R\} \quad \mbox{and} \quad  SO(1,1)= \{ e^{M_{-1} t}, t \in \R\}. $$
A surface of $\C^2$ which is $SO(2)$ or $SO(1,1)$-equivariant may be locally parametrized by an immersion of the form
$$ f(s,t)=e^{M_{\eps_1} t}(z_1(s) , z_2(s) ).$$
We first compute the first derivatives of the immersion:
\begin{eqnarray*} f_s&=&e^{M_{\eps_1} t}(z'_1 , z'_2 ) ,\\
f_t&=&e^{M_{\eps_1} t} M_{\eps_1} (z_1,z_2)=e^{M_{\eps_1} t} (-\eps_1 z_2,z_1).
\end{eqnarray*}
Therefore the Lagrangian condition yields:
\begin{eqnarray*}
0&=&\omega_p(f_s,f_t)= \omega_p((z'_1 , z'_2 ),(-\eps z_2 , z_1 ))\\
&=&- \Im(z'_1 \bar{z}_2)+  \Im(z_2' \bar{z}_1) \\
&=&\frac{d}{ds}\Im(z_2 \bar{z}_1).
\end{eqnarray*}
Hence $z_1 \bar{z}_2$ must be constant. Observe that there is no loss of generality in assuming that $\Im(z_1 \bar{z}_2)$ vanishes: otherwise,
we introduce $\tilde{f}:=\left( \begin{array}{cc} 1 & 0 \\ 0 & e^{i \arg z_2(0)}  \end{array} \right)f,$ which is congruent to $f$.
Thus,  $z_1$ and $z_2$ have the same argument. Next introduce polar coordinates $z_1= r_1 e^{i\varphi}$ and $z_2=r_2 e^{i \varphi}$ and consider separately the definite and indefinite cases:

\medskip

 \noindent \textit{The definite case $p=0$}:

\noindent   The second coordinate of $f$ is
    $$z_2(s) \cos t + z_1(s) \sin t=(r_2(s) \cos t + r_1(s) \sin t)e^{i\varphi(s)}.$$
    Clearly, $\forall s \in I,$ there exists $t(s) \in \R$ such that $r_2(s) \sin t(s) + r_1(s) \cos t(s)=0,$ hence the second coordinate
    of $f$ vanishes at $(s,t(s)).$ Setting $\ga(s):=z_1(s) \cos t(s) - z_2(s) \sin t(s),$ i.e.\ $\ga(s)$ is the first coordinate of $f$ at
    $(s,t(s)),$ we see that  ${f}(s,t)=e^{M_1 (t-s(t))}(\ga(s) ,0).$ Hence the immersion
    $\tilde{f}(s,t):=e^{M_1 t}(\ga(s) ,0)=(\ga \cos t , \ga \sin t)$ parameterizes the same surface as $f$, and we get the required parameterization for the surface $\L.$

  \medskip

  \noindent \textit{The indefinite case $p=1$}:

 \noindent  We first observe that  $r_1 \neq r_2$ since otherwise the immersion would be degenerate.
    If $r_1 > r_2,$  there exists $t(s)$ such that $r_2(s) \cosh t(s) + r_1(s) \sinh t(s)=0,$ hence the second coordinate
    of $f$ vanishes at $(s,t(s)).$ Analogously to the definite case, we set $\ga(s)=r_1(s) \cosh t(s) + r_2(s) \sinh t(s)$, and
    as before, we check that $\tilde{f}(s,t):=(\ga(s) \cosh t, \ga(s) \sinh t)$ parametrizes the same surface as $f.$
    The argument is similar if $r_1 < r_2$: we find $t(s)$ in order to make the first coordinate vanish and find
    $\tilde{f}(s,t)=(\ga(s) \sinh t, \ga(s) \cosh t).$

\medskip

\noindent \textbf{Second case: $n\geq 3$.}
First, set
three different indexes $j,k$ and $l$  and consider the two matrices
$M_{jl}$ and $M_{kl}$ defined by
$$  M_{jl}e_j=e_l, \quad \quad M_{jl}e_l=\eps_j \eps_l e_j \quad \mbox{and} \quad M_{jl}e_m=0 \mbox{ for } m \neq j,l,$$
and
$$  M_{kl}e_k=e_l, \quad \quad M_{kl}e_l=\eps_k \eps_l e_k \quad \mbox{and} \quad M_{kl}e_m=0 \mbox{ for } m \neq k,l.$$
The reader may check that $M_{jl}$ and $M_{kl}$ are skew  with
respect to $\<.,.\>_p.$ Hence, given a point $z$ in $\L,$ the two
curves $s \mapsto e^{M_{jl}s}$ and $s \mapsto e^{M_{kl}s}$ belong
to $SO(p,n-p).$ By the equivariance assumption, it follows that
the curves $s \mapsto e^{M_{jl}s}z$ and $s \mapsto e^{M_{kl}s}z$
belong to $\L$, so the two vectors $M_{jl}z$ and $M_{kl}z$ are
tangent to $\L$ at $z.$ Moreover the Lagrangian assumption yields
$$ 0=\omega_p(M_{jl} z, M_{kl}z)= \Re z_j \Im z_k - \Re z_k \Im z_j.$$
Since this holds for any pair of indexes $(j,k)$, it follows that
$\Re z$ and $\Im z$ are collinear. Therefore there exist $\varphi \in
\R$ and $y \in \R^n$ such that $z= e^{i \varphi} y.$ Let $r
>0$ and $x \in \X^{n-1}_{p,\eps}$ such that $y=rx$, and set $\ga:=r e^{i\varphi}.$ By the
equivariance assumption, the $(n-1)$-dimensional quadric $\ga \,
\X^{n-1}_{p,\eps}$ of $\C^n$ is contained in $\L.$ Finally, since
$\L$ is $n$-dimensional, it must be locally foliated by a
one-parameter family of quadrics $\ga(s)\X^{n-1}_{p,\eps},$ which
proves the first part of the theorem (characterization of
equivariant Lagrangian submanifolds).


We now prove the second part of the theorem: let $f$ be an immersion as described in the statement of the theorem,
  $x$ a point of $\X_{p,\eps}^{n-1}$
and $(e_1, ... , e_{n-1})$  an oriented orthonormal basis of
 $T_x\X_{p,\eps}^{n-1}.$ Setting
 $$  X_j:=\ga e_j  \quad \mbox{ and }\quad X_n:=\ga' x,$$
 it is easy to check that $( X_1, ... , X_{n})$ is a
  basis of $T_{\ga x} \L.$ Then, we calculate
 $$ \omega_p(X_j,X_k)=\<\j X_j,X_k \>_{2p}= \<i\ga ,\ga\> \<e_j,e_k\>_p=0,$$
 $$ \omega_p(X_j,X_n) =\<\j X_j,X_n \>_{2p}= \<i\ga ,\ga' \> \<e_j,x\>_p=0,$$
which shows that
$\L$ is Lagrangian.  Finally, we get the Lagrangian angle of $\L$
as follows:
\begin{eqnarray*}  e^{i \be} &=& \Omega (X_1, ... , X_{n}) \\
&=& \Omega ( \ga e_1, ... , \ga e_{n-1},\ga' x) \\
&=& \ga' {\ga}^{n-1} \Omega ( e_1, ... ,  e_{n-1},x)=\ga' {\ga}^{n-1}.
\end{eqnarray*}

\medskip

From Theorem \ref{equivLagCnTheo}, it is straightforward to
describe equivariant minimal Lagrangian submanifolds: $\be$
vanishes if and only if $\Im \ga' \ga^{n-1}=0$, which we easily integrate to  get $  \Im \ga^n = c$
for some real constant $c.$ If $c$ vanishes, the curve $\ga$ is made up of $n$ straight lines
passing through the origin, and the corresponding Lagrangian submanifold is nothing but
the union of $n$ linear spaces of $\C^n.$ If $c$ does not vanish, the curve is
a made up of $2n$ pieces, each of one contained in an angular sector $\{\varphi_0 < \arg \ga <\varphi_0 + \frac{\pi}{n}\}.$
 If $n=2,$ they are  hyperbolae. Summing up, we have obtained the
following  characterization of equivariant,
minimal Lagrangian submanifolds:
\begin{coro} \label{equivLagCnCoro}Let $\L$ be a connected, minimal Lagrangian submanifold of $(\C^n,\<\<.,.\>\>_p)$ which
is $SO(p,n-p)$-equivariant. Then $\L$ is congruent to an open subset of either an affine Lagrangian $n$-plane, or of the \em Lagrangian catenoid\index{Lagrangian catenoid} \em
$$\left\{   \ga. x \in \C^n \, \big| \,  x \in \X^{n-1}_{p,\eps}, \ga \in \C, \Im \ga^n = c \right\},$$
where $c$ is a non-vanishing real constant.
 \end{coro}


\subsection{Lagrangian submanifolds from evolving quadrics} \label{secquad}


This section describes a class of Lagrangian submanifolds which
generalize the former ones and follows ideas from \cite{Jo}
(see also \cite{LW},\cite{JLT}). Consider a
 a real, invertible $n\times n$ matrix $M$ which is self-adjoint with respect to $\<.,.\>_p,$
i.e.\ $\<Mx,y\>_p=\<x,My\>_p, \, \forall x,y \in \R^n.$

\begin{theo} \label{theoquad} Let $c \in \R$ such that the quadric
$$\s := \left\{x \in \R^n | \,  \<x,Mx\>_p=c   \right\}$$
is a non-degenerate hypersurface of $(\R^n,\<.,.\>_p)$ and $r(s)$ a positive function on an interval $I$ of $\R.$
 Then the immersion
 $$\begin{array}{clcl} f: & I \times \s &  \to & \C^{n} \\
                    & (s,x) & \mapsto &
  r(s)e^{iM s}x
      \end{array}$$
is a Lagrangian  and its Lagrangian angle is given by
$$\beta = \mbox{\emph tr} Ms + \arg \left( c\frac{r' }{r}+i  |Mx|^2_p \right) +\pi/2.$$

\end{theo}


\noindent \textit{Proof.}
 Let $(e_1,...,e_{n-1})$ be an orthonormal basis of
 $ T_x \s = (Mx)^\perp,$
that we complete by $e_n$ in such a way that $(e_1,...,e_n)$ is an oriented, orthonormal basis of $\R^n.$ Hence $e_n$
is collinear to $Mx$ and, setting $\eps_n:=|e_n|_p^2$  we have
$$ Mx = \eps_n \<Mx,e_n\>_p e_n.$$
We obtain a basis of tangent vectors to $f(I \times \s)$ at a
point $z=re^{iMs}x$,
 setting
$$Z_j=e^{iMs} e_j \quad \mbox{and} \quad Z_n =(r'+riM) e^{iMs}x.$$
Using the fact that $e^{iMs} \in U(p,n-p)$,
 it is easily  checked that $\omega_p(Z_j,Z_n) $ and $\omega_p(Z_j,Z_k) $ vanish,
hence the immersion $f$ is Lagrangian. To complete the proof, we compute
 \begin{eqnarray*} \Omega(Z_1,...,Z_n)&=&\Omega(e^{iMs} e_1, ... , e^{iMs} e_{n-1},(r'+irM) e^{iMs}x)\\
&=& i \det_{\C} [e^{iMs}] \det_{\C}(e_1, ..., e_{n-1},(r'+irM)x) \\
&=&i \det_{\C} [e^{iMs}]  \left( r' \eps_n \<x,e_{n}\>_p+ i r \eps_n \<Mx,e_n\>_p \right).
\end{eqnarray*}
Using the fact that
$$\<x,e_n\>_p=\frac{\<x,Mx\>_p}{\eps_n \<Mx,e_n\>_p}=\frac{c}{\eps_n \<Mx,e_n\>_p} ,$$
 we get
$$\Omega(Z_1,...,Z_n)=i e^{i \, \mbox{\rm tr} Ms}  \frac{ r}{ \<Mx,e_n\>_p} \left( c\frac{r'}{r} + i \eps_n \<Mx,e_n\>_p^2  \right).$$
We deduce, using the fact that $\eps_n \<Mx,e_n\>_p^2 = |Mx|_p^2,$
\begin{eqnarray*} \be &=& \arg(\Omega(Z_1,...,Z_n))\\
&=& \frac{\pi}{2}+ \mbox{\rm tr} Ms + \arg \left( c\frac{r'}{r} + i |Mx|_p^2  \right),
\end{eqnarray*}
which is the required formula. 

\begin{exem} \em \label{EquivEx} Assume that $M=Id$ and $c=1.$ Then $f$ becomes
$$ \begin{array}{lccc}  f :
 & I \times \X^{n-1}_{p,\eps}
 &\to&  \C^n \\
& (s,x)& \mapsto &  r(s)e^{is} x.
\end{array}$$
In particular the image of the immersion is a $SO(p,n-p)$-equivariant submanifold as in Section \ref{equivLagCn}.
\em
\end{exem}


\begin{coro} \label{coroquad} The Lagrangian immersion $f$ introduced in  Theorem \ref{theoquad} above is minimal if and only if one of the three statements holds:
\begin{itemize}
    \item[\emph{(i)}] $\mbox{\emph{tr}}\,  M=0$ and the function $r$ is constant;
    \item[\emph{(ii)}] $\mbox{\rm tr} \,  M=0$ and the constant $c$ vanishes;
    \item[\emph{(iii)}] the image of $f$ is a part of the Lagrangian catenoid described in the previous section.
\end{itemize}

\end{coro}

\noindent \textit{Proof.}
The Lagrangian angle $\be$ must be constant, so  the term  $\arg \left( c\frac{r' }{r}+i  |Mx|^2_p \right)$
must be independent of $x.$ This happens if and only if either $r'$ or $c$ vanishes, or both $|Mx|^2_p$ and $\frac{r'}{r}$
are constant. If $r'$ or $c$ vanish, the first term $\mbox{\rm tr} Ms$ of $\be$ must be constant as well, hence we must
have $\mbox{\rm tr}\,  M=0.$ These are the first two cases of the corollary. Suppose now $|Mx|^2_p$ is constant on $\s$, i.e.\
$$ \forall x \in \R^n \mbox{ such that } \<Mx,x\>_p=c,   \, |Mx|_p^2 = c'.$$
Since $M$ is invertible, it is equivalent to
$$ \forall y \in \R^n \mbox{ such that } \<y,M^{-1}y\>_p=c,   \, |y|_p^2 = c'.$$
It follows that the  quadric $\{ \<y,M^{-1}y\>_p=c\}$ is contained in the quadric $\X^{n-1}_{p,c'}$, hence
$M^{-1}$ is a multiple of the identity and so is $M.$ Hence
 the immersion is equivariant and we are in the situation described in Example~\ref{EquivEx} above. The result follows
 from Corollary
 \ref{equivLagCnCoro}.


\begin{exem} \em \label{quadind} Set $n=2$, $p=1$ and $M=\left( \begin{array}{cc} 0 & -1 \\ 1 & 0  \end{array} \right).$
Since $\rm{tr} \, M =0$ the immersion $f(s,x)=r(s)e^{iM s}x$ is minimal if $c$ vanishes or if $r$ is constant.
The case of vanishing $c$ is trivial: the quadric $\s$ reduces to the union of the two straight lines $\{x_1=0\}$ and $\{x_2=0\},$
and the image of $f$ is the union the two complex planes $\{ z_1 = 0 \}$ and $\{ z_2 = 0\}.$

 In the case of non-vanishing $c$, constant $r$, the set
$$\s=\{ x \in \R^2 | \, \<x,Mx\>_1 =2x_1 x_2=c\}$$
is an hyperbola which  may be parametrized by $t \mapsto (e^t, \frac{2}{c} e^{-t}).$ On the other hand
 $$ e^{iMs}=\left( \begin{array}{cc} \cosh s & -i\sinh s \\ i\sinh s & \cosh s  \end{array} \right).$$
 so, setting $r = 1$,   we are left with the immersion
$$f(s,t)=\big(e^t \cosh s - i \frac{2}{c} e^{-t} \sinh s,\frac{2}{c} e^{-t} \cosh s + i e^t \sinh s \big).$$
Observing  that
 $(s,t)$ are conformal coordinates, we obtain null coordinates setting $u:=s+t$ and $v:=s-t$. 
It follows that the immersion takes the form $f(u,v)=\ga_1(u)+ \j \ga_2(v)$ where
$$\ga_1(u):=\frac{1}{2}\big(e^{u}+ \frac{2}{c}i e^{-u},e^{-u}+\frac{2}{c}ie^{u} \big)$$
and
 $$\ga_2(v):= \frac{-1}{2}\big(e^{v}+ i\frac{2}{c}e^{-v},e^{v}+i \frac{2}{c}e^{-v} \big)$$
 are two hyperbolae in the null plane $P=\{x_1-y_2=0, x_2-y_1 =0\}.$
We therefore recover a special case of
 Theorem \ref{minlagn=2}. \em
\end{exem} 

\begin{exem} \em In the definite case, since the metric $\<.,.\>_0$ is positive, there exists an orthonormal basis of eigenvectors of $M$. So we may assume
without loss of generality that $M=diag(\la_1, ...,\la_n)$, where the $\la_j$s are real constants. It follows that a point $z$ of $\L$ takes the form
$$ (x_1 e^{i\la_1 s}, ... , x_n e^{i\la_n s}),$$
where $\sum_{j=1}^n \la_j x_j^2 = c.$
We observe furthermore that $\L$ is a properly immersed submanifold if and only if all the coefficients $\la$ are rationally related. In this case we may assume without loss of generality that they are integer numbers. This case is studied \cite{LW}. Observe moreover that $c$ cannot vanish (otherwise $\s$ reduces to the origin), so by Theorem \ref{coroquad}, $\L$ is minimal if and only if $\rm{tr} \, M =0.$
Example \ref{quadind} above proves that the situation is richer in the indefinite case, since there may not exist an orthonormal basis of eigenvectors.
\em
\end{exem}

\bigskip

\noindent
Henri Anciaux \\
Universidade de S\~ao Paulo \\
  IME, Bloco A \\
  1010 Rua do Mat\~ao  \\
Cidade Universit\'aria   \\
 05508-090 S\~ao Paulo, Brazil   \\
henri.anciaux@gmail.com \\

\end{document}